\documentclass[letterpaper, 10 pt, conference]{ieeeconf}
\IEEEoverridecommandlockouts

\usepackage[utf8]{inputenc}
\usepackage{amsmath}
\usepackage{mathrsfs}
\usepackage{amssymb}
\usepackage{color}
\usepackage{float}
\usepackage{graphicx}
\usepackage{tabularx}
\usepackage{array}
\usepackage[table]{xcolor}
\usepackage{graphicx} 
\usepackage{bernsteinStyle}
\usepackage[colorlinks=true, linkcolor=blue, citecolor=blue, urlcolor=blue]{hyperref}
\usepackage[dvipsnames]{xcolor}

\usepackage[
  backend=biber,
  style=ieee,
  citestyle=ieee-comp,
  sorting=none
]{biblatex}
\DeclareFieldFormat{labelnumberwidth}{#1}
\usepackage{csquotes}
\AtBeginBibliography{\scriptsize\setlength{\itemsep}{0pt}\setlength{\parskip}{0pt}}

\usepackage{xcolor}

\usepackage{tikz}
\definecolor{fblue}{rgb}{.8,.9,1}
\usetikzlibrary{shapes,arrows.meta,calc,positioning,bending}
\usetikzlibrary{calc,patterns,decorations.pathmorphing,decorations.markings,shapes,arrows,positioning,calc,shapes.multipart,backgrounds,fit,quotes}
\usetikzlibrary{plotmarks,chains,calc,shapes.geometric}
\tikzstyle{bigblock} = [draw, fill=fblue, rectangle, 
    minimum height=6em, minimum width=8em]
\tikzstyle{medblock} = [draw, fill=fblue, rectangle, 
    minimum height=4em, minimum width=4em]    
\tikzstyle{mux} = [draw, fill=black!20, rectangle, 
    minimum height=5em, minimum width=0.1em]    
\tikzstyle{smallblock} = [draw, fill=fblue, rectangle, 
    minimum height=3em, minimum width=4em]
\tikzstyle{sum} = [draw, fill=fblue, circle, node distance=1cm]
\tikzstyle{signal} = [coordinate]
\tikzstyle{pinstyle} = [pin edge={to-,thin,black}]
\tikzstyle{block} = [draw, fill=fblue, rectangle, 
    minimum height=3em, minimum width=6em]
\tikzstyle{blockS} = [draw, fill=fblue, rectangle, 
    minimum height=3em, minimum width=4em]    
\tikzstyle{input} = [coordinate]
\tikzstyle{output} = [coordinate]

\tikzset{fontscale/.style = {font=\relsize{#1}}
    }
    
\usepackage{longtable,tabularx}

\title{A Numerical Investigation\\ of Extremum-Seeking-Based Command Generation\\ for Adaptively Controlled Systems}
\author{
Jhon Manuel Portella Delgado,
Aidan Rice,
Jacob C. Vander Schaaf,
and
Dennis S. Bernstein.
\thanks{Jhon Manuel Portella Delgado, Aidan Rice,
Jacob Vander Schaaf,
and Dennis S. Bernstein are with the Department of Aerospace Engineering, University of Michigan, Ann Arbor, MI, USA. {\tt\small \{jhonp,aidancr,jacobcvs, dsbaero\}@umich.edu}}
}
\date{February 2026}

\begin{document}

\maketitle

\begin{abstract}
We develop an adaptive feedback control technique that combines an extremum-seeking-based command generator (ECG) with indirect adaptive control.
In particular, ECG is used to generate commands that asymptotically optimize a cost function that is measured but whose functional form is unknown.
For feedback control with command following and stabilization, the present paper combines ECG with predictive cost adaptive control (PCAC), which is an indirect adaptive control extension of model predictive control (MPC).
PCAC extends generalized predictive control (GPC) 
by using quadratic programming to enforce output constraints and recursive least squares (RLS) with variable-rate forgetting (VRF) for system identification.
The resulting ECG/PCAC framework combines command generation with closed-loop system identification and online optimization.
%
%
The contribution of this paper is a numerical investigation of ECG/PCAC for adaptive stabilization, command following, and disturbance rejection.
\end{abstract}

\textit{\bf keywords:} adaptive control, extremum seeking, predictive cost adaptive control.



\section{Introduction}

The underlying strength of feedback control is the ability to achieve performance despite uncertainty \cite{bernsteinessay}, which may arise due to unknown dynamics and disturbances.
To take advantage of this ability, robust control uses a prior characterization of uncertainty to guarantee performance over the range of uncertainty.
In contrast, adaptive control adjusts the controller to optimize performance for the actual plant and disturbances.
Consequently, adaptive control can provide better performance than robust control while accommodating a larger range of physical plant variations.

In many control applications, the objective is to minimize the command-following error, where a command signal is provided by the user.
In some applications, however, it is more natural to optimize a cost function rather than to specify a command.
The cost function may represent a performance metric that may be measured in real time, but whose functional form is unknown, and thus the optimal command is also unknown.
Relevant problems of technological interest include energy efficiency and lift enhancement \cite{stoutScitech2026}.

Online optimization problems with an unknown performance metric are amenable to extremum seeking \cite{KrsticBook}.
The stability and convergence properties of extremum seeking have been extensively analyzed  \cite{krstic2000,krstic2000_2,tan2010extremum}.
Applications include robotics \cite{matveev2015,bagheri2018,calli2018,sotiropoulos2019}, energy systems \cite{ghaffari2014,ye2016,bizon2017,zhou2018}, and aerospace  \cite{pokhrel2024extremum,pokhrel2024novel,yuan2024extremum}.
The essential feature of extremum seeking is a  dither signal, which enables gradient estimation and determines a search direction toward a local extremum.
In practical implementation, however, the dither signal leads to persistent oscillations, which may be undesirable.
To address this issue, various modifications of extremum seeking have been proposed, including approaches where the amplitude of the dither signal decreases as the system approaches the optimizer or is attenuated by additional dynamics \cite{wang2016,haring2016,haring2018,yin2019,bhattacharjee2021}.
In \cite{paredes2024}, a dither-reduction mechanism is introduced that does not depend on the system state. 
This technique is limited to output minimization and exhibits reduced convergence speed.

The present paper focuses on applications of extremum seeking where the plant is stabilized by a feedback controller and where extremum seeking is used to generate a command signal to be followed by the feedback controller.
The {\it extremum-seeking-based command generator} (ECG) generates a command signal that optimizes a measured but otherwise unknown performance metric.
Variations of ECG have been considered in \cite{krstic2014extremum,raafat2018multivariable}.


For feedback control, this paper combines ECG with predictive cost adaptive control (PCAC) \cite{islamPCAC}, which is an indirect adaptive control method based on model predictive control (MPC) \cite{kwon2006receding,camacho2013model,cairano2018MPC,eren2017model}.
PCAC is an extension of generalized predictive control (GPC) \cite{ClarkeGPC1,ClarkeGPC2,Bitmead:90,Mosca:95}, which uses MPC with closed-loop system identification performed by recursive least squares (RLS) with variable-rate forgetting (VRF).
PCAC extends GPC by using quadratic programming to enforce output constraints and windowed VRF to accelerate learning within RLS.

Unlike data-driven MPC approaches \cite{allgowerlearning2022}, PCAC does not rely on offline data or prior models.
Instead, PCAC performs closed-loop system identification using RLS with VRF \cite{islam2019recursive,adamRLS2020,ankitRLS2020,bruceNandSRLS,mohseni2022recursive}.
PCAC relies on self-generated excitation for system identification \cite{islamPCAC} and therefore it does not require probing signals as in dual control \cite{mesbahMPCdual2}.
In addition, as in MPC, PCAC enforces input and output constraints through quadratic programming, including magnitude and rate limits.
Applications of PCAC include noise and vibration control \cite{rezaMSSP,PCACRijkeACC2025,rileyBACT2024}, atmospheric flight control \cite{rileyQCACC2025,rileyScitech2025}, and flow control \cite{PCACAirfoilACC2024}, with experimental results in \cite{rezaMSSP,PCACRijkeACC2025,rileyQCACC2025}.


%
%
The resulting ECG/PCAC framework combines command generation with closed-loop system identification and online optimization.
ECG/PCAC thus enables simultaneous stabilization and optimal command generation without requiring prior models or offline training data.
The contribution of this paper is a numerical investigation of the synergy between ECG for performance optimization and PCAC for adaptive stabilization, command following, and disturbance rejection.

\section{Problem Formulation}
\label{sec:Prob_Formulation}
\begin{figure}[!ht]
\centering
\resizebox{\columnwidth}{!}{%
\begin{tikzpicture}[>={stealth'}, line width=0.25mm, node distance=2em]

\tikzstyle{block} = [draw, rounded corners, minimum height=3em, minimum width=7em, align=center]
\tikzstyle{smallblock} = [draw, rounded corners, minimum height=3em, minimum width=5em, align=center]
\tikzstyle{sum} = [draw, circle, inner sep=0.3em]

\node[draw=none] at (0,0) (orig) {};

\node[draw, fill=green!20, rounded corners, minimum width=12cm, minimum height=6cm] 
at ([yshift=0.25em, xshift=4.8em]orig.center) (M) {};

\node[block, fill=blue!20] (mpc) at ([xshift=-3cm, yshift=0.8cm]M.center)
{
MPC
};

\node[block, fill=yellow!20, below=2.5em of mpc] (id)
{
RLS/VRF
};

\node[smallblock, fill=blue!10, right=3.5em of mpc] (zoh)
{
ZOH
};

\node[block, fill=red!20, right=3.5em of zoh] (sys)
{
Plant
};

\node[
    draw,
    rounded corners,
    fill=red!40,
    fill opacity=0.25,
    text opacity=1,
    fit=(mpc)(id),
    inner xsep=0.35cm,
    inner ysep=0.5cm
] (pcacbox) {};

\node[anchor=north west] at ([xshift=0.08cm,yshift=-0.05cm]pcacbox.north west) {\small PCAC};

\node[sum,minimum size=0.5cm] (sumrk) at ([xshift=-1cm]mpc.west) {};

\node at ([xshift=-0.42cm,yshift=0.22cm]sumrk.center) {\scriptsize $+$};
\node at ([xshift=0.2cm,yshift=-0.32cm]sumrk.center) {\scriptsize $-$};

\coordinate (rin) at ([xshift=-0.6cm]sumrk.west);
\draw[->] (rin) -- node[above, xshift=-0.8cm] {$r_k$} (sumrk.west);

\draw[->] (sumrk.east) -- node[above,xshift=-0.2 cm] {$e_k$} (mpc.west);

\draw[->] (mpc.east) -- node[above] {$u_k$} (zoh.west);

\draw[->] (zoh.east) -- node[above] {$u$} (sys.west);

\coordinate (fbdown) at ([xshift=1.75cm,yshift=-1.4cm]sys.south);
\coordinate (fbleft) at ([xshift=3.5cm]id.east);
\coordinate (swL)    at ([xshift=-1.10cm]fbleft);
\coordinate (swR)    at ([xshift=-0.35cm]fbleft);

\draw[-] (sys.east)--([xshift = 0.5cm]sys.east) -- node[right, pos=0.35] {$y$} (fbdown);
\draw[-] (fbdown) -- (swR);
\draw[->] (swL) -- node[above] {$y_k$} (id.east);

\draw[fill=white] (swL) circle (0.045cm);
\draw[fill=white] (swR) circle (0.045cm);

\draw[-] (swL) -- ([yshift=0.18cm]swR);

\node[above=0.07cm of $(swL)!0.5!(swR)$] {$T_\rms$};

\coordinate (sumdown) at ([yshift=-1.4cm]sumrk.south);
\draw[-] ([xshift=-0.7cm,yshift=-1.5cm]swL) -- ([xshift=-0.7cm]swL);
\draw[-] ([xshift=-0.7cm,yshift=-1.5cm]swL) -- ([xshift=0cm,yshift=-1.8cm]sumdown);
\draw[->] ([xshift=0cm,yshift=-1.8cm]sumdown) -- (sumrk.south);

\coordinate (uTap)  at ([xshift=0.9cm]mpc.east);
\coordinate (uDown) at ([yshift=-1.65cm]uTap);

\draw[-] (uTap) -- (uDown);
\draw[->] (uDown) -- ([yshift=0.3cm]id.east);

\draw[->] (id.north) -- node[right] {$\theta_k$} (mpc.south);

\coordinate (Jout) at ([yshift = 0.2cm,xshift=1.5cm]sys.east);
\draw[-] ([yshift = 0.2cm]sys.east) -- node[above, xshift=0.5cm] {$J$} (Jout);
\draw[->]([yshift = 0.8 cm]sys.north) -- node[above, xshift=0.2cm] {$w$} (sys.north);

\node[block, fill=cyan!20, minimum width=6em] (esc) at ([yshift=-4.1cm]M.center)
{
ECG
};

\coordinate (swECG_R) at ([xshift=2.3cm]esc.east);
\coordinate (swECG_L) at ([xshift=1.5cm]esc.east);

\coordinate (outerRight) at ([xshift=0.0cm]Jout);
\coordinate (outerDown)  at ([yshift=-5.1cm]outerRight);

\draw[-] (Jout) -- (outerRight);
\draw[-] (outerRight) -- (outerDown);
\draw[-] (outerDown) -- (swECG_R);

\draw[fill=white] (swECG_L) circle (0.045cm);
\draw[fill=white] (swECG_R) circle (0.045cm);

\draw[-] (swECG_L) -- ([yshift=0.18cm]swECG_R);

\node[above=0.07cm of $(swECG_L)!0.5!(swECG_R)$] {$T_\rms$};

\draw[->] (swECG_L) -- node[above] {$J_k$} (esc.east);

\coordinate (outerLeft)  at ([xshift=-5.7cm]esc.west);
\coordinate (outerUp)    at ([xshift=-0.7cm]rin);

\draw[-] (esc.west) -- (outerLeft);
\draw[-] (outerLeft) -- (outerUp);
\draw[-] (outerUp) -- (rin);

\end{tikzpicture}
}
\caption{ECG/PCAC sampled-data control architecture for adaptive stabilization, command following, and disturbance rejection.}
\label{fig:PCAC_ESC}
\end{figure}
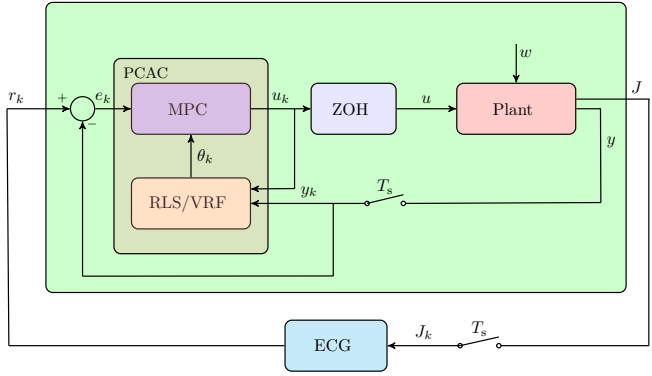

We consider the closed-loop control architecture shown in Figure \ref{fig:PCAC_ESC}, where the goal is to optimize the measured output $J \in \BBR^{p}$ of the plant.
Using sampled values $J_k$ of $J,$ ECG determines the command $r_k \in\BBR^{m}$ for the feedback control system. 
Since the functional form of $J$ is unknown, analytical optimization of $J$ is not possible.

\section{Extremum-Seeking-Based Command Generator}
\label{sec:ECG_controller}

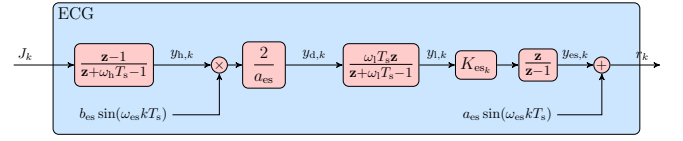
\begin{figure}[!ht]
    \centering
    \resizebox{\columnwidth}{!}{%
    \begin{tikzpicture}[>={stealth'}, line width = 0.25mm]
    \node[draw = none] at (0,0) (orig) {};
    \node [smallblock, rounded corners, minimum height = 3.25cm, minimum width = 14.5cm] at ([yshift = -0.25em, xshift = 5.8em]orig.center) (esc_controller) {};
    \node[below right] at (esc_controller.north west) {\large$\rm ECG$};
    \node[smallblock, fill=red!20, rounded corners, minimum height = 3em, minimum width = 3em] 
    at (orig.center) (LP) 
    {\large$\dfrac{2}{a_{\rm es}}$};
    \node[sum, fill=red!20, inner sep = 0.4em, left = 1em of LP.west] (mult_c) {};
    \node[draw = none] at (mult_c.center) {$\times$};
    \node[smallblock, fill=red!20, rounded corners, minimum height = 3em, minimum width = 3em, left = 4em of mult_c.west] (HP) {\Large$\frac{\bfz - 1}{\bfz + \omega_{\rmh}T_\rms - 1}$};
    \node[smallblock, fill=red!20, rounded corners, minimum height = 3em, minimum width = 3em, right = 4em of LP.east] (ECG_int) {\Large$\frac{\omega_\rml T_\rms\bfz}{\bfz + \omega_\rml T_\rms - 1}$};
    \node[smallblock, fill=red!20, rounded corners, minimum height = 2em, minimum width = 2.5em, right = 2.5em of ECG_int.east] (scale2a) {\large$K_{{\rm{es}}_k}$};
    %
    \node[smallblock, fill=red!20, rounded corners, minimum height = 2em, minimum width = 2em, right = 1.5em of scale2a.east] (Kesc) {\Large$\frac{\bfz}{\bfz - 1}$};
    \node[sum, fill=red!20, inner sep = 0.4em, right = 2.5em of Kesc.east] (sum_c) {};
    \node[draw = none] at (sum_c.center) {$+$};
    \node[draw = none] at ([xshift = -6.75em, yshift = -3.5em]mult_c.center) (ditherECG) {$b_{\rm es} \sin(\omega_{\rm es} k T_\rms)$};
    \node[draw = none] at ([xshift = 27em]ditherECG.center) (ditherECG_2) {$a_{\rm es} \sin(\omega_{\rm es} k T_\rms)$};
    %
    %
    \draw[->] (ditherECG_2.east) -| (sum_c.south);
    \draw[->] (ditherECG.east) -| (mult_c.south);
    \draw[->] (HP.east) -- node[above] {$y_{\rmh,k}$} (mult_c.west);
    \draw[->] (mult_c.east) -- (LP.west);
    \draw[->] (LP.east) -- node[above] {$y_{\rmd,k}$} (ECG_int.west);
    \draw[->] (ECG_int.east) -- node[above] {$y_{\rml,k}$}(scale2a.west);
    \draw[->] (scale2a.east) -- (Kesc.west);
    \draw[->] (Kesc.east) -- node[above] {$y_{{\rm es},k}$} (sum_c.west);
    \draw[->] ([xshift = -4.25em]HP.west) -- node[above, xshift = -1.3em] {$J_k$}(HP.west);
    \draw[->] (sum_c.east) -- node[above, xshift = 0.5em] {$r_k$}([xshift = 3.5em]sum_c.east);
    \end{tikzpicture}
    }
    \caption{
    Extremum-seeking-based command generator ECG with measured performance $J_k$ and command $r_k \in \BBR.$
    }
    \label{fig:ESC}
\end{figure}

Figure \ref{fig:ESC} shows that ECG is based on the extremum-seeking controller described in \cite{paredes2024retrospective,portella2026extremum}.
The dynamics of ECG are given by
\begin{align}
    y_{{\rm{h}},{k}}
        &=
            (1-
                \omega_{\rm{h}}
                T_\rms
            )
            y_{{\rm{h}},{k-1}}
            +
            J_k
            -
            J_{k-1},
    \\
    y_{{\rmd},k}
        &=
            \tfrac
            {
            2b_{\rm{es}}
            }
            {
            a_{\rm{es}}
            }
             \omega_{\rml}
            T_\rms
            \sin
            (
                k\omega_{\rm{es}}
                T_\rms
            )
             y_{\rmh,k},
    \\
    y_{{\rm{l}},k}
        &=
            (
                1-\omega_{\rml}
                T_\rms
            )
            y_{\rml,{k-1}}
            +
            \omega_{\rml}
            T_{\rms}y_{\rmd,{k}},
    \\
    y_{{\rm{es}},k}
        &=
            y_{{\rm{es}},{k-1}}
            +
            K_{{\rm es},{k}}
            \,
            y_{{\rml},k},
    \\
    r_k
        &=
            y_{{\rm{es}},k}
            +
            a_{\rm{es}}
            \sin
            (
                k\omega_{\rm{es}}
                T_\rms
            ),
    \label{eq:ECG_output}
\end{align}
where $y_{{\rm h},k}$ is the output of the highpass filter at step $k$, $y_{{\rm l},k}$ is the output of the lowpass filter at step $k$, and $y_{{\rmd},k}$ is an internal variable coming out from the demodulation process.
Note that the approximation of the gradient of $J$ at step $k$ is $y_{{\rm l},k}.$
Moreover, $y_{\rm es}$ denotes the output of the integrator and $K_{{\rm es},k}$ is the gradient gain, while $\omega_{\rm l} > 0$ and $\omega_{\rm h} > 0$ are the cutoff frequencies of the lowpass and highpass filters, respectively. 
The parameter 
$\omega_{\rm es} > 0$ is the dither signal frequency, $a_{\rm es}$ is the modulation dither signal amplitude, and $b_{\rm es}$ is the demodulation dither signal amplitude.

Note that $K_{\rm es}$ is positive when the objective is to maximize $J$ and negative when the objective is to minimize $J$. 
The signal $a_{\rm es} \sin(\omega_{\rm es} t)$ corresponds to the \textit{modulation dither signal}, while $b_{\rm es} \sin(\omega_{\rm es} t)$ corresponds to the \textit{demodulation dither signal}. 
Finally, the highpass filter removes the bias from $y$, and the lowpass filter
removes high-frequency components from $y_{{\rmd},k},$ to produce $y_{{\rml},k},$ which approximates the gradient \cite{krstic2000,KrsticBook2003}.
In this work, $y_{\rmh,0} = y_{\rml,0} = y_{{\rm{es}},0} = 0,$ and we define $e_k \isdef r_k - y_k.$

\section{Predictive Cost Adaptive Control} \label{sec:PCAC}
Predictive cost adaptive control (PCAC) combines online identification  with output-feedback MPC. 
The PCAC algorithm is presented in this section.
Subsection \ref{subsec:ID} describes the technique used for online identification, namely, recursive least squares (RLS). 
Subsection \ref{subsec:MPC} reviews the output-feedback MPC technique for receding-horizon optimization.
%

\subsection{Online Identification Using Recursive Least Squares} \label{subsec:ID}

Let $\hat n\ge 0$ and, for all $k\ge 0,$ let $\hat{F}_{1,k},\ldots,\hat{F}_{\hat n,k} \in\BBR^{p\times p}$ and $\hat{G}_{0,k},\ldots,\hat{G}_{\hat n,k}\in\BBR^{p\times m}$  be the coefficient matrices to be estimated using RLS.
Furthermore, let $\hat y_k\in\BBR^{p}$ be an estimate of the output $y_k \in \BBR^{p}$ defined  by
\begin{align}
    \hat{y}_{k} = -\sum_{i=1}^{\hat n}{\hat F}_{i,k} y_{k-i} + \sum_{i=0}^{\hat n}{\hat G}_{i,k} u_{k-i}, \label{eq:yhat}
\end{align}
where  
\begin{gather}
   y_{-\hat n}=\cdots= y_{-1}=0,\\ u_{-\hat n}=\cdots=u_{-1}=u_{0}=0,
\end{gather}
and $u_k \in \BBR^{m}$.
The input-output model \eqref{eq:yhat} can be written as
\begin{align}
    \hat y_k = \theta_k \phi_k,
\end{align}
where the matrix of unknown parameters $\theta_k \in \BBR^{ p \times [ \hat n( p + m ) + m ]},$ and the regressor vector $\phi_k \in \BBR^{ \hat n ( p + m ) + m },$ are defined as
\begin{align} 
    \theta_k 
        &\isdef 
            \matl 
                -\hat F_{1,k} 
                & 
                \cdots 
                & 
                -\hat F_{\hat n,k} 
                & 
                \hat G_{0,k} 
                & 
                \cdots 
                & 
                \hat G_{\hat n,k} 
            \matr,
    \\
    \phi_k 
        &\isdef 
            \matl 
                y_{k-1} 
                & 
                \cdots 
                & 
                y_{k-\hat n} 
                & 
                u_{k} 
                & 
                \cdots 
                & 
                u_{k-\hat n}  
            \matr^\rmT.
\end{align}
%
%
\begin{align}
    z_k(\theta_k) = y_k  - \hat y_k = y_k  - \theta_k \phi_k
\end{align}

For all $j = 1,\ldots,p,$ we denote each row of $\theta_k$ as $\theta_{j,k} \in \BBR^{ 1 \times [ \hat n ( p + m ) + m ] }$.

For online identification, RLS is used to estimate the coefficients of the input-output model \eqref{eq:yhat}.
To do this, RLS minimizes the cumulative cost 
\begin{align}
    J_k(\hat\theta) 
        = 
            \sum_{j=1}^{p}  J_{j,k}(\hat\theta_j),
    \label{eq:Jprecursor}
\end{align}
where
\begin{align}
    J_{j,k}(\hat\theta_j) 
        \isdef 
            \sum_{i=0}^k 
            \frac{\rho_i}{\rho_k}z_{j,i}(\hat\theta_j)^2
            + 
            \frac{1}{\rho_k}(\hat\theta_j-\theta_{j,0})P_0^{-1}
            (\hat\theta_j-\theta_{j,0} )^\rmT,
    \label{eq:J}
\end{align}
where, for all 
$
    k\ge0,$  $\rho_k \isdef \prod_{j=0}^k \lambda_j^{-1} \in \mathbb{R},
$ 
$
    \lambda_k\in(0,1]
$ 
is the forgetting factor,  
$
    P_0\in\BBR^{[\hat{n}(m + p) + m] \times [\hat{n}(m + p) + p]}
$ 
is positive definite, and
$
    \hat\theta 
        \isdef 
            \matl 
                \hat\theta_1^\rmT 
                & 
                \cdots 
                & 
                \hat\theta_{p}^\rmT 
            \matr^\rmT
            \in\BBR^{p \times [\hat{n}(m + p) + m]}.
$
For all 
$j = 1,\ldots,p,$ 
$
    \theta_{j,0}\in\BBR^{1 \times [\hat{n}(m + p) + m]}
$ 
is the initial estimate of $j$-th row of the coefficient matrix, and
the output-wise performance variable $z_{j,k}(\theta_{j,k} ) \in\BBR$ is defined by
\begin{align}
     z_{j,k}(\theta_{j,k}) 
        \isdef 
            y_{j,k} 
            - 
            \theta_{j,k}\phi_k. 
    \label{eq:zid1}
\end{align}
Note that, with \eqref{eq:zid1}, the cost function \eqref{eq:Jprecursor} is strictly convex and quadratic, and thus has a unique global minimizer.
Assuming that $P_0 = \overline P_0 I_p$, where $\overline P_0 \in \BBR$ is a tuning parameter, the unique global minimizer is computed by RLS using 
\begin{align} 
    L_{k} 
        &= 
            \lambda_{k}^{-1} P_{k}, 
    \label{eq:RLS1} 
    \\
    P_{k+1} 
        &=  
            L_k  
            - 
            \dfrac
            {
                1
            }
            {
                1 
                + 
                \phi_k^\rmT   
                L_k 
                \phi_k
            }  
            L_k 
            \phi_k  
            \phi_k^\rmT L_k, 
    \label{eq:RLS2} 
    \\
    \theta_{k+1} 
        &= 
            \theta_{k}
            +
            \left(
                y_{k}  
                -   
                \theta_{k}
                \phi_k  
            \right)
            \phi_k^\rmT 
            P_{k+1}. 
    \label{eq:RLS3}
\end{align}
Note that $\theta_{k+1}$ computed using  \eqref{eq:RLS3} is available at step $k,$ and thus, $\hat{F}_{1,k+1},\ldots,\hat{F}_{\hat n,k+1}, \hat{G}_{0,k+1},\ldots,\hat{G}_{\hat n,k+1}$ are available at step $k$.
%

In the present paper, the exponential-resetting variant of RLS (ER-RLS) is used.
As discussed in \cite{brianResettingRLS}, ER-RLS modifies the RLS update scheme such that information matrix $R_k = P_k^{-1}$ converges to a specified information matrix $R_\infty$ in the absence of persistent excitation.
This prevents $P_k$ from exhibiting covariance windup.
Due to a required inversion of $R_k$, this algorithm runs in $\SO (n^3)$ where $n$ is the number of estimated parameters.
By comparison, RLS runs in $\SO (pn^2)$ where $n$ is the number of estimated parameters and $p$ is the number of outputs.

Reconsider RLS with the modified information matrix update
\begin{align}
    R_{k+1}  = \lambda_k R_k + (1-\lambda_k) R_\infty + \phi_k \phi_k^\rmT
\end{align}
where positive definite $R_\infty\in\BBR^{[\hat{n}(m+p)+m]\times [\hat{n}(m+p)+m]}$ is a tunable parameter, and retaining the $\theta_k$ update
\begin{align} 
    \theta_{k+1} &= \theta_{k}+    ( y_{k}  -   \theta_{k}  \phi_k  )\phi_k^\rmT R_{k+1}^{-1}.
\end{align} 
As shown in \cite{brianResettingRLS}, this scheme ensures that, if $R_0$ is positive definite, then for all $k \geq 0$, $R_k$ is positive definite.
Additionally, a tight lower bound for $R_k$, dependent on step $k$, exists, and under no excitation, the sequence $\{ R_k\}_{k=0}^\infty$ converges to $R_\infty$.
If $\SO (pn^2)$ complexity must be maintained, Cyclic-Resetting RLS (CR-RLS) \cite{brianResettingRLS} may instead be used, however in this case only the subsequence of $P_k$ every $n$ steps converges to $P_\infty$, and sequence of the covariance matrix oscillates but remains bounded close to $P_\infty$ under zero excitation.
%

%

The step-dependent parameter $\lambda_k$ is the {\it forgetting factor.}
In the case where $\lambda_k$ is constant, RLS uses {\it constant-rate forgetting} (CRF);  otherwise, RLS uses {\it variable-rate forgetting} (VRF).
For F-Test VRF \cite{mohseni2022recursive}, $\lambda_k$ is given by
\begin{align}
    \lambda_k 
        = 
            \dfrac
            {
                1
            }
            {
                1
                +
                \eta 
                g(z_{k-\tau_{\text{d}}
            }, 
            \dots, z_{k}) 
            {\bf 1}
            [g(z_{k-\tau_\text{d}}, \dots, z_{k})]},
    \label{eq:VRF}
\end{align}
where ${\bf 1}\colon \mathbb{R} \rightarrow \{0, 1\}$ is the unit step function, and 
\begin{align}
    &g(z_{k-\tau_\text{d}}, \dots, z_{k})
    \nn \\
        &\isdef  
            \sqrt
            {
                \dfrac
                {
                    \tau_\mathrm{n}
                }
                {
                    \tau_\text{d}
                } 
                \dfrac
                    {
                        \textrm{tr}
                            \left(
                                \Sigma_{\tau_\mathrm{n}}(z_{k-\tau_{\mathrm{n}}}
                                , 
                                \dots
                                , 
                                z_{k})
                                \Sigma_{\tau_\text{d}}(z_{k-\tau_\text{d} }
                                , 
                                \dots
                                , 
                                z_{k})^{-1}
                            \right)
                    }
                    {
                        c
                    }
            }
    \nn\\ 
    & -\sqrt
        {
            F_{p \tau_n,b}
            (1 - \alpha)
        },
\end{align}
where $\eta > 0$ and $p \leq \tau_\mathrm{n} < \tau_\text{d}$ represent numerator and denominator window lengths. 
$\Sigma_{\tau_{\text{n}}}$ and $\Sigma_{\tau_{\text{d}}}$ are the sample variances of the respective window lengths, and
\begin{align}
    a 
        &\isdef 
            \dfrac
            {
                (\tau_n + \tau_d - p - 1)
                (\tau_d - 1)
            }
            {
                (\tau_d - p - 3)
                (\tau_d - p)
            }, 
    \nn \\
    b 
        &\isdef 
            4 
            + 
            \dfrac
            {
                (p \tau_n + 2)
            }
            {
                (a - 1)
            }, 
    \nn\\
    c 
        &\isdef 
            \dfrac
            {
                p \tau_n (b-2)
            }
            {
                b(\tau_d-p-1)
            }.
\end{align}

\subsection{Model Predictive Control (MPC)} \label{subsec:MPC}
We define the tracking output $y_{\rmt,k}\in \BBR$ as
\begin{align}
    y_{\rmt,k} 
        \isdef 
            C_\rmt 
            y_k.  
    \label{Ct}
\end{align}  
The performance objective is to have  $y_{\rmt,k}\in\BBR^{{p}_\rmt}$ follow a commanded sequence of $r_k\in\BBR^{{p}_\rmt},$ whose future values may or may not be known.
In addition to the performance objective, the constrained output $y_{\rmc,k}\in\BBR^{{p}_\rmc}$ is defined by
\begin{align}
    y_{\rmc,k} 
        \isdef 
            C_\rmc 
            y_k, 
    \label{eq:yc}
\end{align}
where $C_\rmc\in\BBR^{{p}_\rmc\times p}.$
The objective is to enforce the inequality constraint 
\begin{align}
    \SC y_{\rmc,k} 
    + 
    \SD 
        \leq 
            0_{n_\rmc\times 1}, 
    \label{eq:constraint}
\end{align}
where $\SC \in\BBR^{n_\rmc\times {p}_\rmc}$ and $\SD \in\BBR^{n_\rmc}.$
Note that \eqref{eq:constraint}, where ``$\le$'' is interpreted component-wise, defines a convex set. 

The control is constrained in both magnitude and rate.
The magnitude control constraint has the form
\begin{align}
    u_{\rm min} 
        \le 
            u_k 
                \le 
                    u_{\rm max}, 
    \label{eq:magsat}
\end{align}
where $u_{\rm min}\in\BBR^{m}$ is the vector of the minimum control magnitudes  and $u_{\rm max}\in\BBR^{m}$ is the vector of maximum control magnitudes.
In addition, the increment-size (rate) control constraint has the form
\begin{align}
    \Delta u_{\rm min} 
        \le 
            u_{k}-u_{k-1} 
                \le 
                    \Delta u_{\rm max}, 
    \label{eq:movesat}
\end{align}
where $\Delta u_{\rm min}\in\BBR^{m}$ is the vector of minimum control increment sizes and $\Delta u_{\rm max}\in\BBR^{m}$ is the vector of maximum control increment sizes.
Note that \eqref{eq:magsat} and \eqref{eq:movesat}, where ``$\le$'' is interpreted component-wise, define convex sets.
  
Next, let $\ell\ge 1$ be the horizon and, for all $k\ge0$ and all $i=1,\ldots,\ell,$ let  $\hat{y}_{k|i}\in\BBR^{p}$ be the $i$-step predicted output, and $u_{k|i}\in\BBR^{m}$ be the $i$-step predicted control.
Then, the $\ell$-step predicted output of \eqref{eq:yhat} for a sequence of $\ell$ future controls is given by
\begin{align}  
    Y_{k,\ell} 
        = 
            \Gamma_k 
            + 
            T_k 
            U_{k,\ell}. 
    \label{ypredict}
\end{align}
where
\begin{align}
        Y_{k,\ell} 
            &\isdef 
                \matl 
                    y_{ k|1} 
                    & 
                    \cdots 
                    & 
                    y_{k|\ell}  
                \matr^\rmT 
        \in \BBR^{\ell p}, 
        \nn \\
        U_{k,\ell} 
            &\isdef 
                \matl 
                    u_{k|1} 
                    & 
                    \cdots 
                    & 
                    u_{k|\ell}  
                \matr^\rmT 
        \in \BBR^{\ell m},
        \nn \\ 
        D_{\hat n,k} 
            &\isdef 
                \matl    
                    y_{k-\hat n + 1} 
                    & 
                    \cdots 
                    & 
                    y_k 
                    & 
                    u_{k-\hat n + 1} 
                    & 
                    \cdots 
                    & 
                    u_k 
                \matr^\rmT
        \nn \\
        & \quad
        \in \BBR^{ \hat n ( p + m ) }, 
        \label{Dk}
\end{align}
\begin{align}
    \Gamma_k 
        &\isdef 
            \matl  
                -  F_{\rmp,k}^{-1} F_{\rmd,k}  
                & 
                F_{\rmp,k}^{-1} G_{\rmd,k}  
            \matr 
            D_{\hat n,k} 
        \in \BBR^{\ell p}, 
    \nn\\
    T_k 
        &\isdef 
            F_{\rmp,k}^{-1} G_{\rmp,k}  
    \in \BBR^{\ell p \times \ell m},  
    \label{inverses}
\end{align}
\begin{align}
    F_{\rmd,k} 
        &\isdef 
            \matl 
                \hat F_{\hat n,k} 
                & 
                \cdots   
                &
                \hat F_{1,k} 
                    \\
                \vdots 
                & 
                \ddots   
                &
                \vdots  
                    \\
                0_{p\times p} 
                & 
                \cdots   
                &
                \hat F_{\hat n,k} 
                    \\
                0_{p\times p}   
                & 
                \cdots   
                & 
                0_{p\times p} 
                    \\
                \vdots 
                & 
                \ddots   
                &
                \vdots  
                    \\
                0_{p\times p}      
                & 
                \cdots 
                &
                0_{p\times p}
            \matr 
            \in\BBR^{\ell p \times  \hat n p }, 
    \nn \\
    G_{\rmd,k} 
        &\isdef 
            \matl 
                \hat G_{\hat n,k} 
                & 
                \cdots   
                &
                \hat 
                G_{1,k} 
                    \\
                \vdots 
                & 
                \ddots   
                &
                \vdots  
                    \\
                0_{p\times m} 
                & 
                \cdots   
                &
                \hat G_{\hat n,k} 
                    \\
                0_{p\times m}   
                & 
                \cdots   
                & 
                0_{p\times m} 
                    \\
                \vdots 
                & 
                \ddots   
                &
                \vdots  
                    \\
                0_{p\times m}      
                & 
                \cdots 
                &
                0_{p\times m}
            \matr 
            \in\BBR^{\ell p \times  \hat n m }, \label{FdGd} 
\end{align}
\begin{align}
    F_{\rmp,k} 
        &\isdef
            \matl
                I_p       
                & 
                \cdots   
                & 
                0_{p\times p} 
                & 
                0_{p\times p} 
                & 
                \cdots   
                & 
                0_{p\times p}   
                    \\
                \vdots       
                & 
                \ddots   
                & 
                \vdots  
                & 
                \vdots  
                & 
                \ddots   
                & 
                \vdots  
                    \\
                \hat F_{\hat n -1 , k}       
                & 
                \cdots  
                &
                I_p 
                & 
                0_{p\times p}   
                & 
                \cdots   
                & 
                0_{p\times p}  
                    \\
                \hat F_{\hat n,k}    
                & 
                \cdots   
                & 
                \hat F_{1,k}  
                & 
                I_p  
                & 
                \cdots 
                & 
                0_{p\times p}   
                    \\ 
                \vdots       
                & 
                \ddots   
                & 
                \vdots  
                & 
                \vdots 
                & 
                \ddots   
                & 
                \vdots  
                    \\
                0_{p\times p} 
                &
                \cdots 
                &  
                \hat F_{\hat n, k} 
                & 
                \hat F_{\hat n-1, k}   
                &
                \cdots 
                & 
                I_p   
            \matr 
            \nn\\ 
            & \quad
            \in \BBR^{\ell p \times \ell p}, 
            \label{Fp}
\end{align}
\begin{align}
     G_{\rmp,k} 
        &\isdef
            \matl
                \hat G_{0,k}       
                & 
                \cdots   
                & 
                0_{p\times m} 
                & 
                0_{p\times m} 
                & 
                \cdots   
                & 
                0_{p\times m}   
                    \\
                \vdots       
                & 
                \ddots   
                & 
                \vdots  
                & 
                \vdots  
                & 
                \ddots   
                & 
                \vdots  
                    \\
                \hat G_{\hat n -1 , k}       
                & 
                \cdots  
                &
                \hat G_{0,k} 
                & 
                0_{p\times m}   
                & 
                \cdots   
                & 
                0_{p\times m}  
                    \\
                \hat G_{\hat n,k}    
                & 
                \cdots   
                & 
                \hat G_{1,k}  
                & 
                \hat G_{0,k}  
                & 
                \cdots  
                & 
                0_{p\times m}   
                    \\ 
                \vdots       
                & 
                \ddots   
                & 
                \vdots  
                & 
                \vdots 
                & 
                \ddots   
                & 
                \vdots  
                    \\
                0_{p\times m} 
                &
                \cdots 
                &  
                \hat G_{\hat n, k} 
                & 
                \hat G_{\hat n-1, k}   
                &
                \cdots 
                & 
                \hat G_{0,k}   
            \matr 
            \nn\\ 
            & \quad
            \in \BBR^{\ell p \times \ell m}. 
            \label{Gp}
\end{align}
To facilitate the weighting of integrated-command-following error, we define the integrated-command-following-error vector
\begin{align}
     I_{\rmt,k,\ell} 
        \isdef 
            \matl 
                \rmi_{k|1}^\rmT 
                & 
                \rmi_{k|2}^\rmT   
                & 
                \cdots  
                & 
                \rmi_{k|\ell-1}^\rmT  
                & \rmi_{k|\ell}^\rmT  
            \matr^\rmT,
\end{align} 
where
\begin{align}
    \rmi_{k|1} 
        \isdef 
            i_k 
            + 
            y_{\rmt,k} 
            - 
            r_k, 
                \quad 
                    i_k 
                        \isdef 
                            \sum_{i=0}^k y_{\rmt,i} - r_i,
\end{align}
and, for all $i=2,\ldots,\ell,$
\begin{align}
    \rmi_{k|i} 
        \isdef 
            i_{k|i-1} 
            + 
            y_{\rmt,k|i-1} 
            - 
            r_{k|i-1}.
\end{align}
Let the vector 
$
    \SR_{k,\ell} 
        \isdef 
            \matl 
                r_{k+1}^{\rm T} 
                & 
                \cdots 
                & 
                r_{k+\ell}^{\rm T} 
            \matr^{\rm T} 
            \in \BBR^{\ell p_\rmt}
$ 
be composed of $\ell$ future commands,
let $y_{\rmt,k|i} \isdef C_\rmt y_{k|i} \in\BBR^{p_\rmt}$ be the $i$-step predicted command-following output, let 
$
    Y_{\rmt , k , \ell } 
        \isdef 
            \matl 
                y_{ \rmt, k|1}^{\rm T} 
                & 
                \cdots 
                & 
                y_{ \rmt, k|\ell}^{\rm T} 
            \matr^{\rm T} 
                = 
                    C_{\rmt,\ell} 
                    Y_{k,\ell}
                    \in \BBR^{\ell p_\rmt},
$ 
where 
$
    C_{\rmt,\ell} 
        \isdef 
            I_\ell 
            \otimes 
            C_\rmt 
                \in \BBR^{\ell p_\rmt \times \ell p},
$
and define
\begin{align}
    \Delta U_{k,\ell}  
        \isdef 
            \matl
                (u_{k|1}-u_k)^\rmT 
                & 
                (u_{k|2}-u_{k|1})^\rmT 
                & 
                \nn\\
                \cdots 
                & 
                (u_{k|\ell}-u_{k|\ell-1 })^\rmT
            \matr^\rmT 
            \in \BBR^{\ell m}.
\end{align}
Then, the receding horizon optimization problem is given by
\begin{align}
    &\min_{U_{k,\ell}} 
        \left(
            Y_{\rmt,k,\ell} 
            - 
            \SR_{k,\ell}
        \right)^{\rm T}  
        Q 
        \left(
            Y_{\rmt,k,\ell} 
            - 
            \SR_{k,\ell}
        \right)  
        \nn \\
    & 
    + 
    I_{\rmt,k,\ell}^{\rm T}  
    Q_\rmi 
    I_{\rmt,k,\ell}^{}
    +  
    U_{k,\ell}^{\rm T} 
    R  
    U_{k,\ell}  
    + 
    \Delta U_{k,\ell}^{\rm T} 
    R_{\delta} 
    \Delta U_{k,\ell}  
    + 
    \varepsilon^\rmT 
    S 
    \varepsilon  
    \label{cost}
\end{align}
subject to
\begin{align}
    \SC_{\ell}Y_{k,\ell} 
    + 
    \SD_{\ell}  
        \leq 
            \varepsilon, 
    \\
    U_{\rm min} 
        \leq 
        U_{k,\ell}  
        \leq 
        U_{\rm max},  
    \\
    \Delta U_{\rm min} 
        \leq 
        \Delta U_{ k,\ell}  
        \leq 
        \Delta U_{\rm max},  
    \\
    0_{\ell n_\rmc \times 1}  
        \leq 
            \varepsilon,  
    \label{lastconstraint}
\end{align}
where 
$Q   \in \BBR^{ \ell p_\rmt \times \ell p_\rmt }$ is the positive-definite command-following error weighting,  
$Q_\rmi   \in \BBR^{ \ell p_\rmt \times \ell p_\rmt }$ is the positive-semidefinite (PSD) integrated command-following error weighting,  
$R \in \BBR^{\ell m \times \ell m}$ is the PSD control weight, 
$R_\delta \in \BBR^{\ell m \times \ell m}$ is the PSD control increment-size weight, 
$S\in \BBR^{\ell n_\rmc \times \ell n_\rmc}$ is the PSD constraint relaxation weight,
$
    U_{\min} 
        \isdef 
            1_{\ell } 
            \otimes 
            u_{\min} 
            \in \BBR^{\ell m},
$ 
$
    U_{\max} 
        \isdef 
            1_{\ell } 
            \otimes 
            u_{\max} 
            \in \BBR^{\ell m},
$ 
$
    \Delta U_{\min} 
        \isdef 
            1_{\ell } 
            \otimes 
            \Delta 
            u_{\min} 
            \in \BBR^{\ell m},
$
$
    \Delta U_{\max} 
        \isdef 
            1_{\ell } 
            \otimes 
            \Delta 
            u_{\max} 
            \in 
            \BBR^{\ell m},
$ 
and 
$
    {\SC}_{\ell} 
        \isdef 
            I_\ell 
            \otimes 
            (\SC C_\rmc) 
            \in\BBR^{\ell n_\rmc \times \ell p}
$ 
and 
$
    \SD_\ell 
        \isdef 
            1_{\ell \times 1} 
            \otimes 
            \SD 
            \in\BBR^{\ell n_\rmc}.
$
The quadratic programming (QP) optimization \eqref{cost}--\eqref{lastconstraint} is solved using the MATLAB quadprog routine with warm starting of the solution at each step using the solution from the previous step.

In summary, at each time step,  online identification is performed to find input-output model coefficients $\theta_{k+1}$, which are then used to create the matrices $\Gamma_k$ and $T_k$.
The matrices $\Gamma_k$ and $T_k$ are then used in a receding horizon optimization problem to solve for the $\ell$-step  controls $U_{k,\ell}.$
The control input for the next step is  given by $u_{k|1}$, and the remaining components of $U_{k,\ell}$ are discarded.



\section{Numerical Investigation of ECG/PCAC}
\label{sec:examples}

In this section, we apply ECG/PCAC to an undamped oscillator, which is Lyapunov stable, a double integrator, which is quadratically unstable, and an exponentially unstable plant.
All three plants are modeled by linear, second-order, continuous-time systems operating under sampled-data control.
%


\subsection{Undamped oscillator}

Consider the undamped oscillator
\begin{align}
    \dot{x}_1
        &=
            x_2,
    \label{eq:USM_x1_dot}
    \\
    \dot{x}_2
        &=
            -\frac{k}{m}x_1
            +
            \frac{1}{m}u,
    \label{eq:USM_x2_dot}
\end{align}
where $x_1$ is the position of the particle, $x_2$ is the velocity, $m$ is the particle mass, $k$ is the spring stiffness, and $u$ is the force control input.
Note that \eqref{eq:USM_x1_dot}, \eqref{eq:USM_x2_dot} is undamped, and thus it is Lyapunov stable but not asymptotically stable.
The measured output is the position $y = x_1.$
%

\begin{figure}[!ht]
    \centering
    \includegraphics[width=\columnwidth]{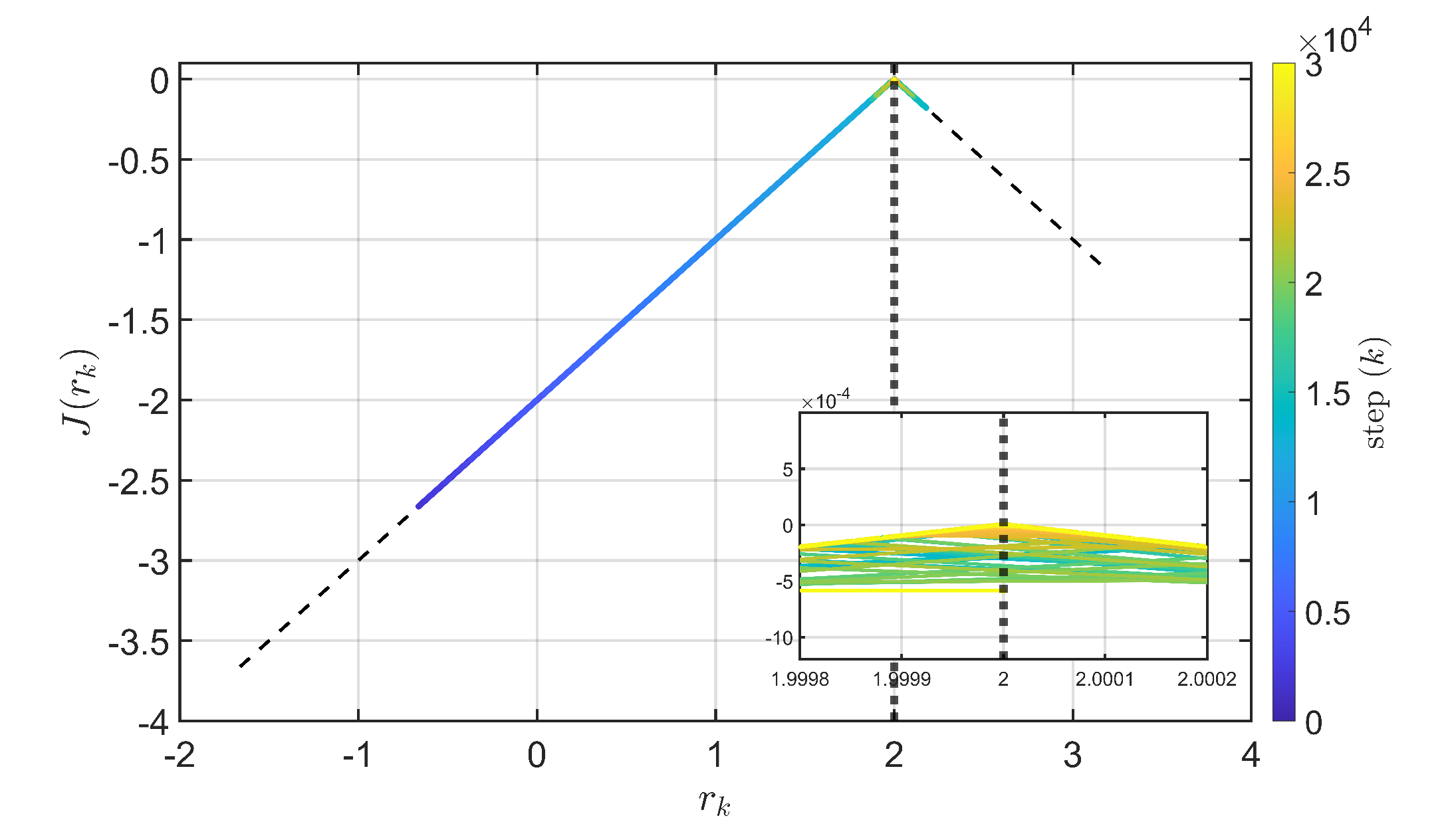}
    \caption{
    \textbf{Undamped oscillator:}
    This plot shows the values of $J(r_k)$ computed using the command $r_k$ obtained from ECG.
    The step-heat map shows that $r_k$ converges to the optimal command $r^\star=2$ denoted by the vertical dotted line.
    %
    The cost function $J(r_k)$ is shown by the dashed black line.
    }
    \label{fig:UND_cost}
\end{figure}

The cost function to be maximized is given by
\begin{align}
    J(r)
        &=
            -|r-r^\star|,
    \label{eq:UD_oscillator_cost_function}
\end{align}
which is maximized by $r=r^\star$, where $J(r^\star)=0$.  
For this example, $r^* = 2.$
Figure \ref{fig:UND_cost} shows that ECG/PCAC updates the command $r$ so that the output of the plant controlled by PCAC converges to $r^*.$  
Although measurements of $J$ are assumed to be available, the functional form of $J$ and the maximizing command $r^*$ are assumed to be unknown.

\begin{figure}[!ht]
    \centering
    \includegraphics[width=\columnwidth]{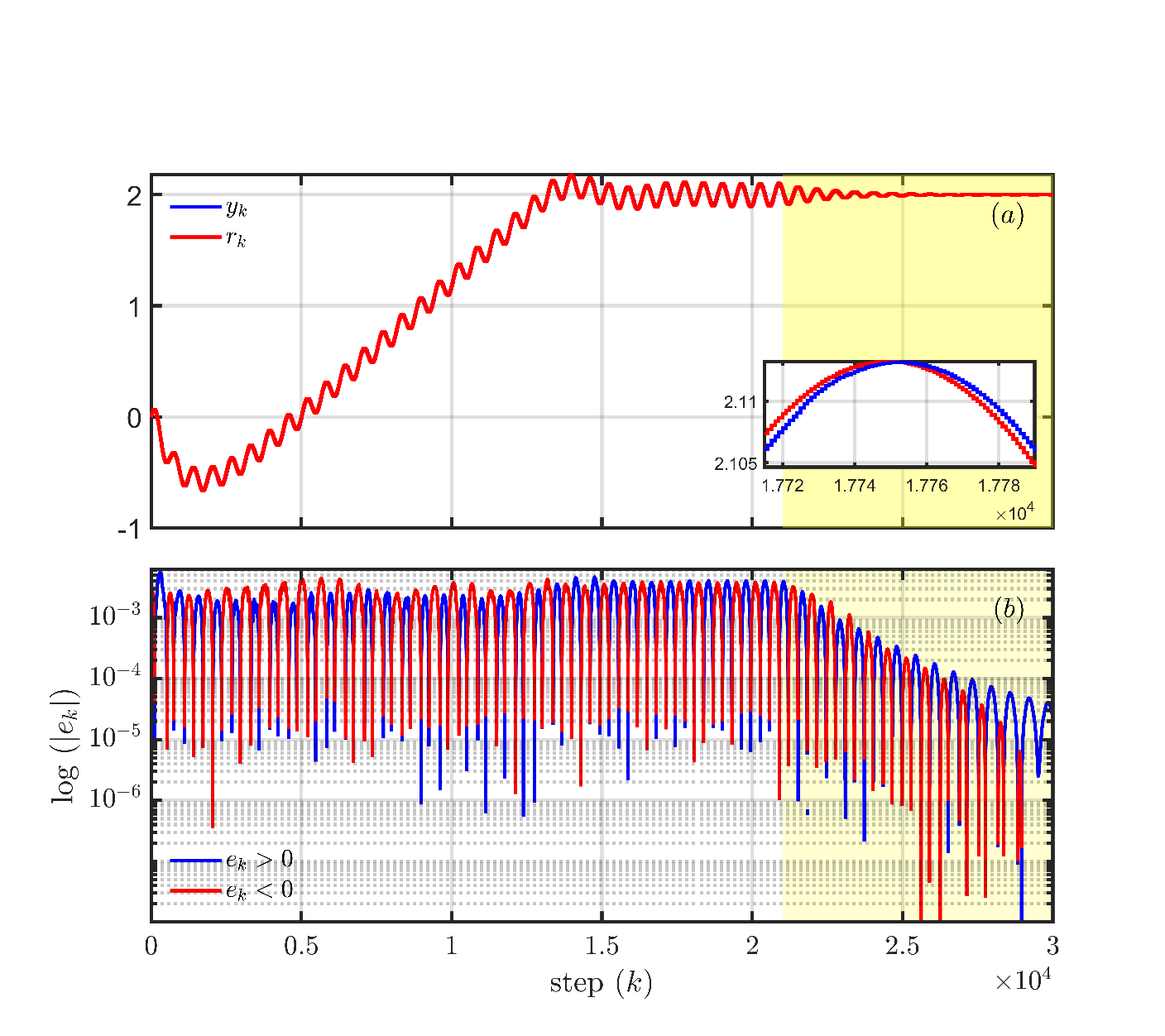}
    \caption{
    \textbf{Undamped oscillator:}
    (a) shows the closed-loop response of \eqref{eq:USM_x1_dot}, \eqref{eq:USM_x2_dot} with $u_k$ given by PCAC and $r_k$ given by ECG.
    The inset and (b) show that PCAC follows $r_k$.
    The dither input from ECG is set to zero at step $k=2.1e4.$
    }
    \label{fig:UND_closed_loop}
\end{figure}

\begin{figure}[!ht]
    \centering
    \includegraphics[width=\columnwidth]{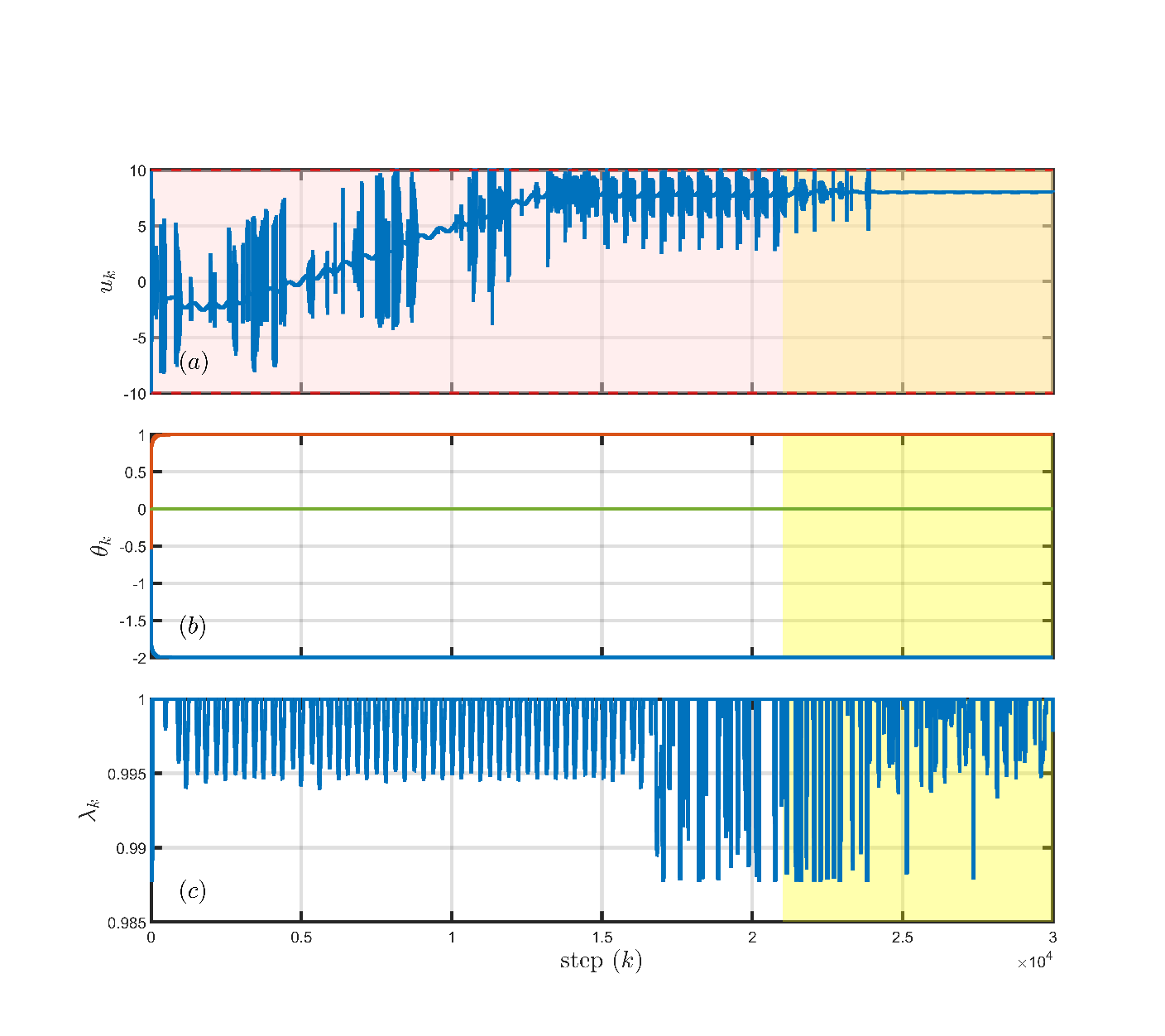}
    \caption{
    \textbf{Undamped oscillator:}
    (a) shows the control signal $u_k$ obtained from PCAC,
    (b) shows the corresponding model coefficients obtained from RLS,
    and (c) shows the variable-rate forgetting factor.
    }
    \label{fig:UND_control}
\end{figure}

ECG/PCAC is implemented with 
sampling time $T_\rms = 10^{-2}$ s, 
dither amplitude $a_{\rm es} = 0.1$, dither  frequency $\omega_{\rm es} = 1\,\mathrm{rad/s}$,
gradient gain $K_{{\rm es},{k}} = 0.1$, and lowpass and highpass cutoff frequencies $\omega_{\rm h} = 0.001\, \mathrm{rad/step}$ and $\omega_{\rm l} = 0.01\, \mathrm{rad/step}$.
Unless noted otherwise,  these values are used for this and all subsequent examples for $k \leq 2.1e4$.
For $k > 2.1e4,$ corresponding to the yellow  regions in Figures \ref{fig:UND_closed_loop} and \ref{fig:UND_control}, the dither amplitude is modulated according to
\begin{align}
    a_{{\rm es},{k+1}} =
\begin{cases}
a_{{\rm es},{k}}, & |y_{{\rm{l}},k}| \geq \beta ,\\
\max \{ a_{{\rm es},\min}, \alpha a_{{\rm es},{k}}\},
& |y_{{\rm{l}},k}| < \beta ,
\end{cases}
%
%
%
    \label{eq:a_es_adaptation}
\end{align}
where
$a_{{\rm es},{k}}$ is the dither amplitude at step $k,$ 
$a_{{\rm es},{\min}} > 0$ is a prescribed lower bound that prevents the dither amplitude from vanishing, $\alpha \in (0,1)$ is a decay factor that reduces the dither amplitude when attenuation is activated, and $\beta  > 0$ is a threshold that determines whether the system is sufficiently close to steady state.
In this example, we set
$
    a_{{\rm es},\min}
        = 5e\mbox{-5}, 
$
$
    \beta 
        = 
            0.02,
$
and
$
    \alpha
        =
            0.999.
$
This dither-amplitude modification is performed after step $k=2.1e4$  to reduce the amplitude of the oscillations due to the dither signal.

PCAC is implemented with model order $\hat{n} = 2$ and horizon $\ell = 30.$   
The state weighting matrix is $Q_i = 0$, and 
the initial covariance is $P_0 = 1e6$. 
The control input is constrained by the magnitude-saturation constraint $u \in [-10,10]$, and no constraints are imposed on the control rate.
Unless noted otherwise,  these values are used for this and all subsequent examples.

Figure \ref{fig:UND_closed_loop} shows the closed-loop response of the undamped oscillator \eqref{eq:USM_x1_dot}, \eqref{eq:USM_x2_dot} under PCAC, together with the command $r_k$ generated by the ECG scheme \eqref{eq:ECG_output} and the command-following error $|e_k|$.
%
%
%
The yellow region shows that the modulation \eqref{eq:a_es_adaptation} asymptotically reduces the oscillations.

Figure \ref{fig:UND_control} shows $u_k$ obtained from the constrained optimization \eqref{cost}--\eqref{lastconstraint}, the estimated model coefficients $\theta_k$ given by \eqref{eq:RLS3}, and the forgetting factor $\lambda_k$ given by \eqref{eq:VRF}.
Figure \ref{fig:UND_cost} shows the values of $J(r_k)$ computed from \eqref{eq:UD_oscillator_cost_function} using the command $r_k$ generated by \eqref{eq:ECG_output}, together with its convergence to $r^\star$.

\subsection{Double integrator}

Consider the double integrator
\begin{align}
    \dot{x}_1
        &=
            x_2,
    \label{eq:DI_x1_dot}
    \\
    \dot{x}_2
        &=
            u + w,
    \label{eq:DI_x2_dot}
\end{align}
where $x_1$ is the position of the particle, $x_2$ is the velocity, $u$ is the force control input, and $w = 0.1$ is a constant disturbance.
Note that \eqref{eq:DI_x1_dot}, \eqref{eq:DI_x2_dot} is quadratically unstable.
The measured output is the position $y = x_1.$

\begin{figure}[!ht]
    \centering
    \includegraphics[width=\columnwidth]{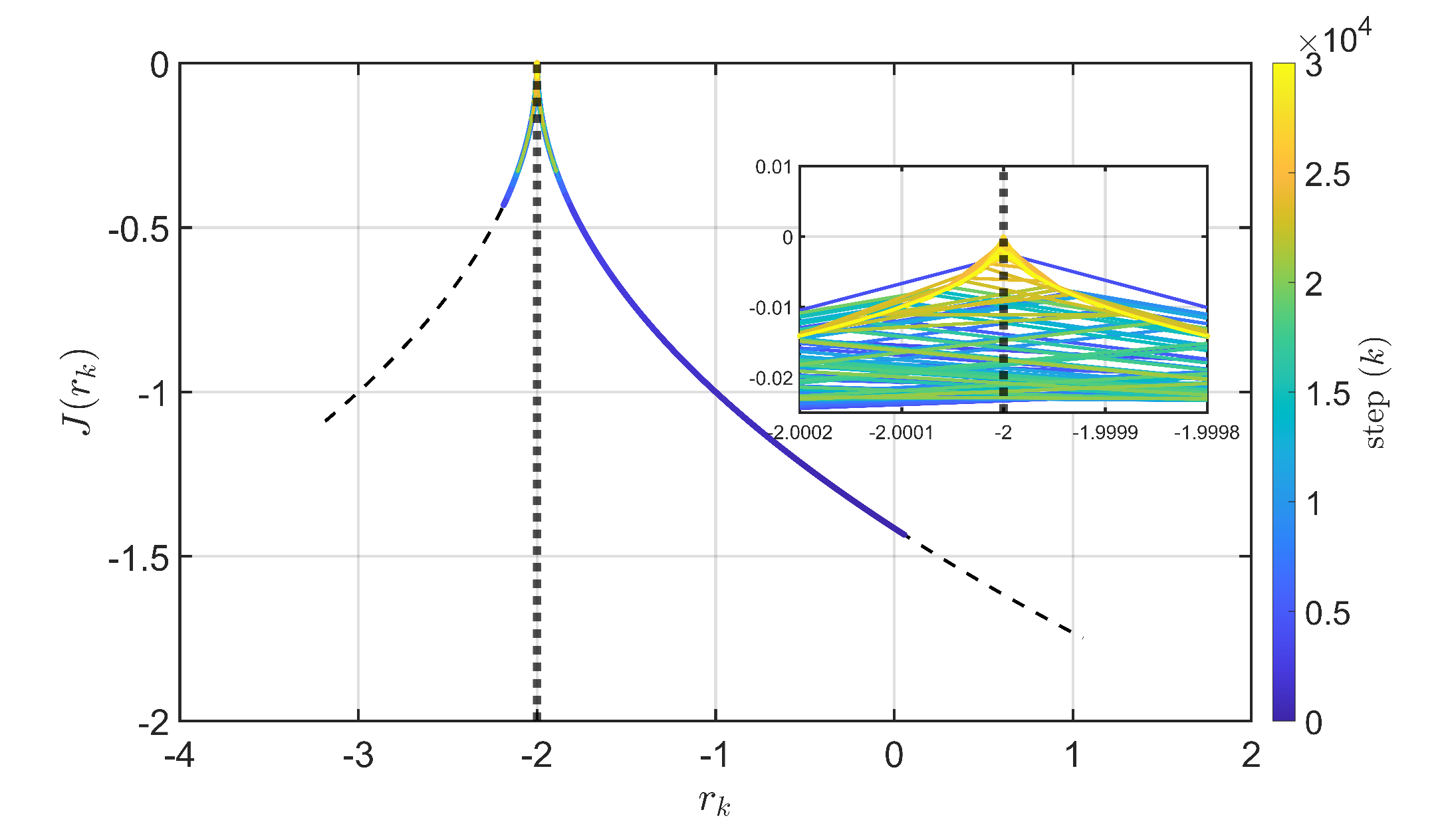}
    \caption{
    \textbf{Double integrator with constant disturbance.} This plot shows the values of $J(r_k)$ computed using the command $r_k$ obtained from \eqref{eq:ECG_output}.
    %
    %
    The step-heat map shows that $r_k$ converges to  $r^\star$ denoted by the vertical dotted line.
    The cost function $J(r_k)$ is shown by the dashed black line.
    }
    \label{fig:Cost_Function_ESC_PCAC_DoubleInt_system}
\end{figure}

The cost function to be maximized is given by
\begin{align}
    J(r)
        &=
            -
            \sqrt
            {
                |r-r^\star|
            },
\label{eq:DI_cost_function}
\end{align}
which is maximized by $r=r^\star$, where $J(r^\star)=0$.  
For this example, $r^* = -2.$
Figure \ref{fig:Cost_Function_ESC_PCAC_DoubleInt_system} shows that ECG/PCAC updates the command $r$ so that the output of the plant controlled by PCAC converges to $r^*.$  
Although measurements of $J$ are assumed to be available, the functional form of $J$ and the maximizing command $r^*$ are assumed to be unknown.

\begin{figure}[!ht]
    \centering
    \includegraphics[width=\columnwidth]{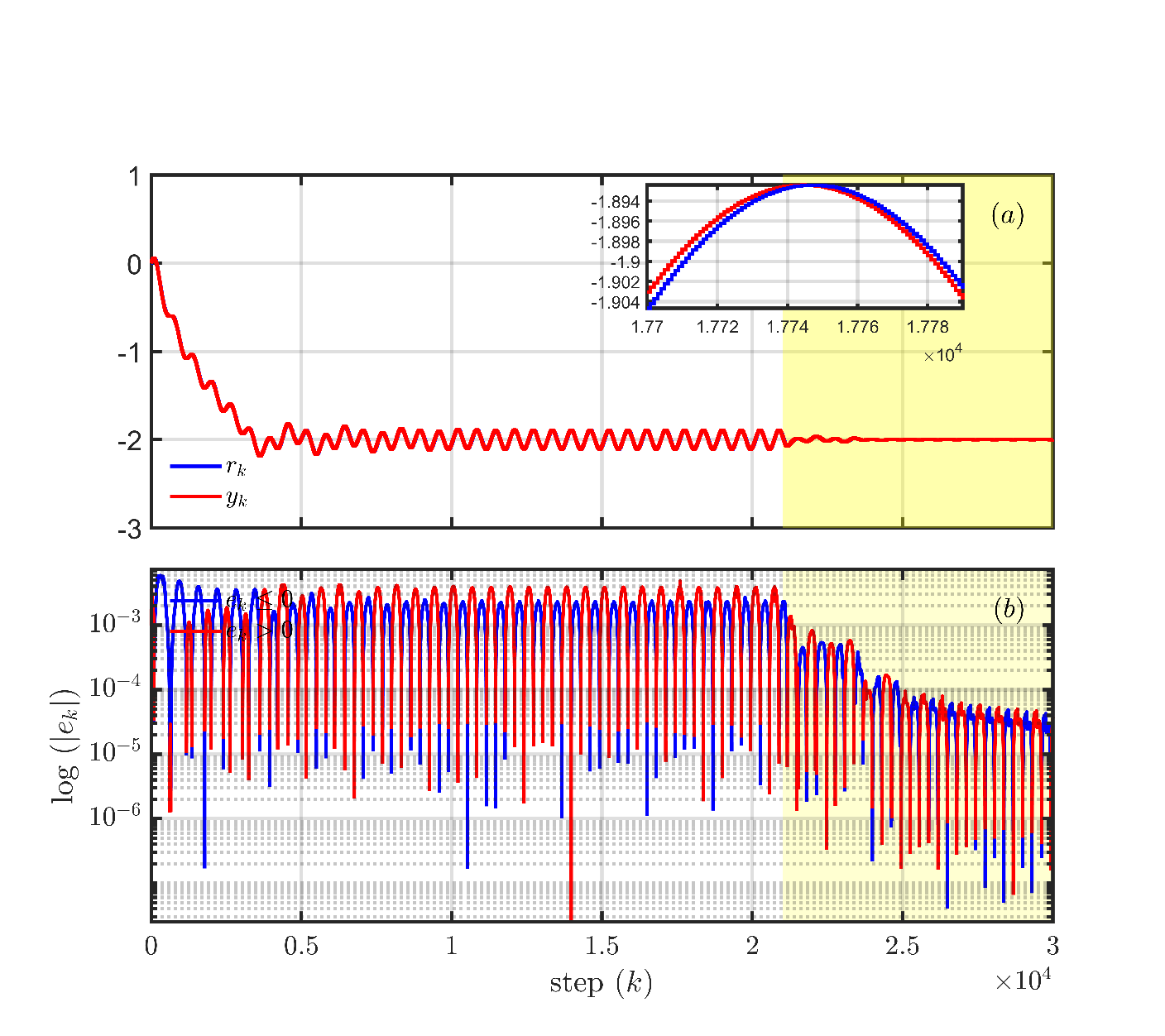}
    \vspace{-.4in}
    \caption{
    \textbf{Double integrator with constant disturbance:} 
    (a) shows the closed-loop response of the double integrator with a constant disturbance \eqref{eq:DI_x1_dot}, \eqref{eq:DI_x2_dot} with the control $u_k$ given by PCAC and the command $r_k$ given by ECG.  
    %
    %
    The inset and (b) show that PCAC follows $r_k$.
    %
    }
    \label{fig:ECG_PCAC_Output_DoubleInt_system}
\end{figure}

\begin{figure}[!ht]
    \centering
    \includegraphics[width=\columnwidth]{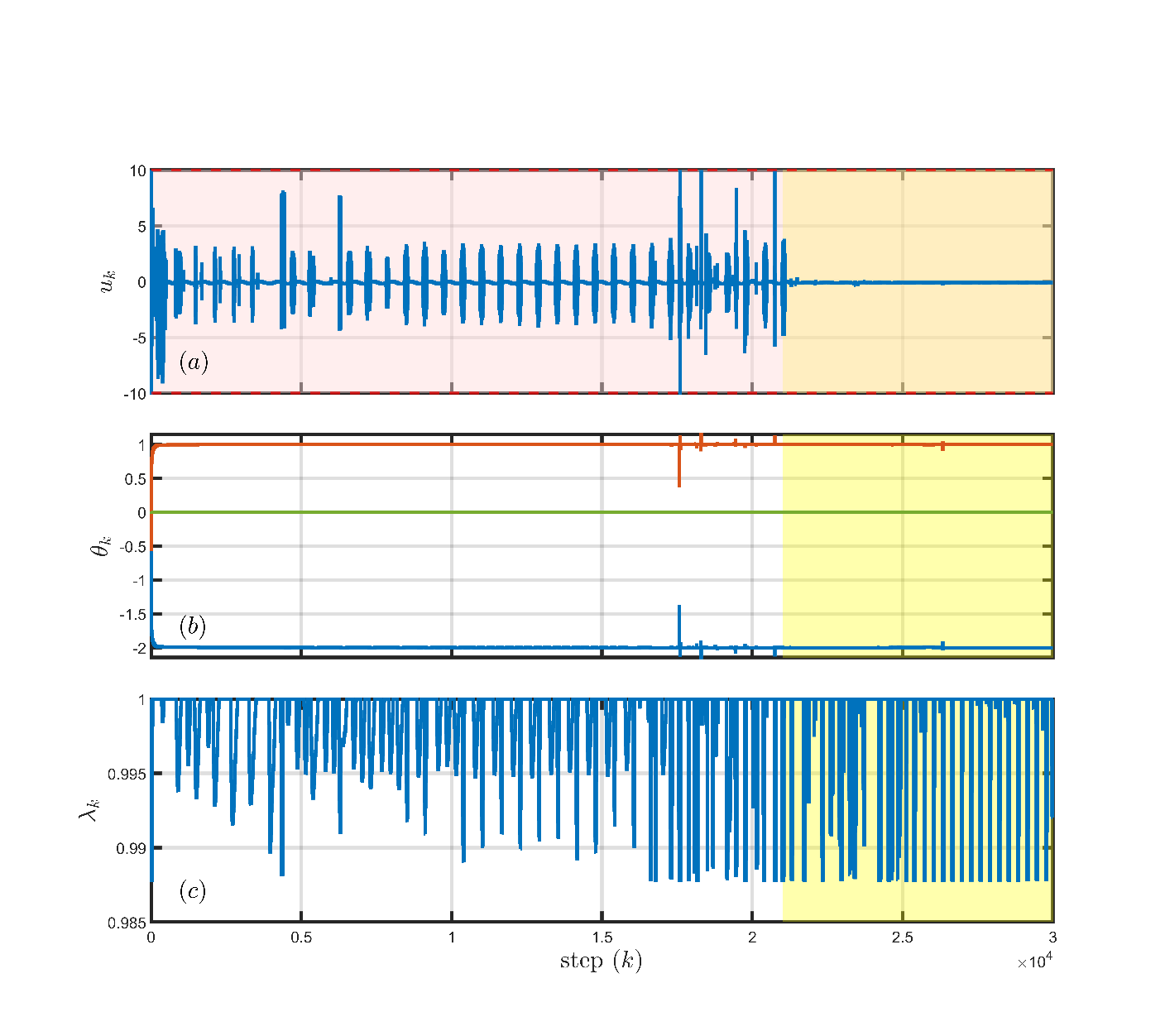}
    \vspace{-.4in}\caption{
    \textbf{Double integrator with constant disturbance.}
    (a) shows the control signal $u_k$ obtained from PCAC,
    (b) shows the corresponding model coefficients obtained from RLS,
    and (c) shows the forgetting factor.
    %
    %
    %
    }
    \label{fig:PCAC_Control_DoubleInt_system}
\end{figure}

For $k > 2.1e4,$ corresponding to the yellow  regions in Figures \ref{fig:ECG_PCAC_Output_DoubleInt_system} and \ref{fig:PCAC_Control_DoubleInt_system}, the dither amplitude is modulated according to \eqref{eq:a_es_adaptation} with 
$
    a_{{\rm es},{\min}} 
        =
            5e\mbox{-4},
$ 
$
    \beta 
        = 
            0.05,
$
and
$
    \alpha = 0.995.
$
This dither-amplitude modification is performed after $2.1e4$ steps to reduce oscillations due to the dither signal.

Figure \ref{fig:ECG_PCAC_Output_DoubleInt_system} shows the closed-loop response of the double integrator \eqref{eq:DI_x1_dot}, \eqref{eq:DI_x2_dot} 
under PCAC, together with the command $r_k$ generated by ECG \eqref{eq:ECG_output} and the command-following error $|e_k|$.
%
%
The yellow region shows that   \eqref{eq:a_es_adaptation}  asymptotically reduces the oscillations.

Figure \ref{fig:PCAC_Control_DoubleInt_system} shows $u_k$ obtained from solving the constrained optimization problem \eqref{cost}--\eqref{lastconstraint}, the estimated model coefficients $\theta_k$ given by \eqref{eq:RLS3}, and the forgetting factor $\lambda_k$ given by \eqref{eq:VRF}.
Figure \ref{fig:Cost_Function_ESC_PCAC_DoubleInt_system} shows the values of $J(r_k)$ computed from \eqref{eq:DI_cost_function} using the command $r_k$ generated by \eqref{eq:ECG_output}, together with its convergence to $r^\star$.

\subsection{Exponentially unstable system}

Consider the exponentially unstable system
\begin{align}
    \dot{x}_1
        &=
            x_2,
    \label{eq:EUS_x1_dot}
    \\
    \dot{x}_2
        &=
            x_1
            +
            u,
    \label{eq:EUS_x2_dot}
\end{align}
where $x_1$ is the position of the particle, $x_2$ is the velocity of the particle, and $u$ is the force control input.
Note that \eqref{eq:EUS_x1_dot}, \eqref{eq:EUS_x2_dot} is exponentially unstable.
The measured output is the position $y = x_1.$
%

\begin{figure}[!ht]
    \centering
    \includegraphics[width=\columnwidth]{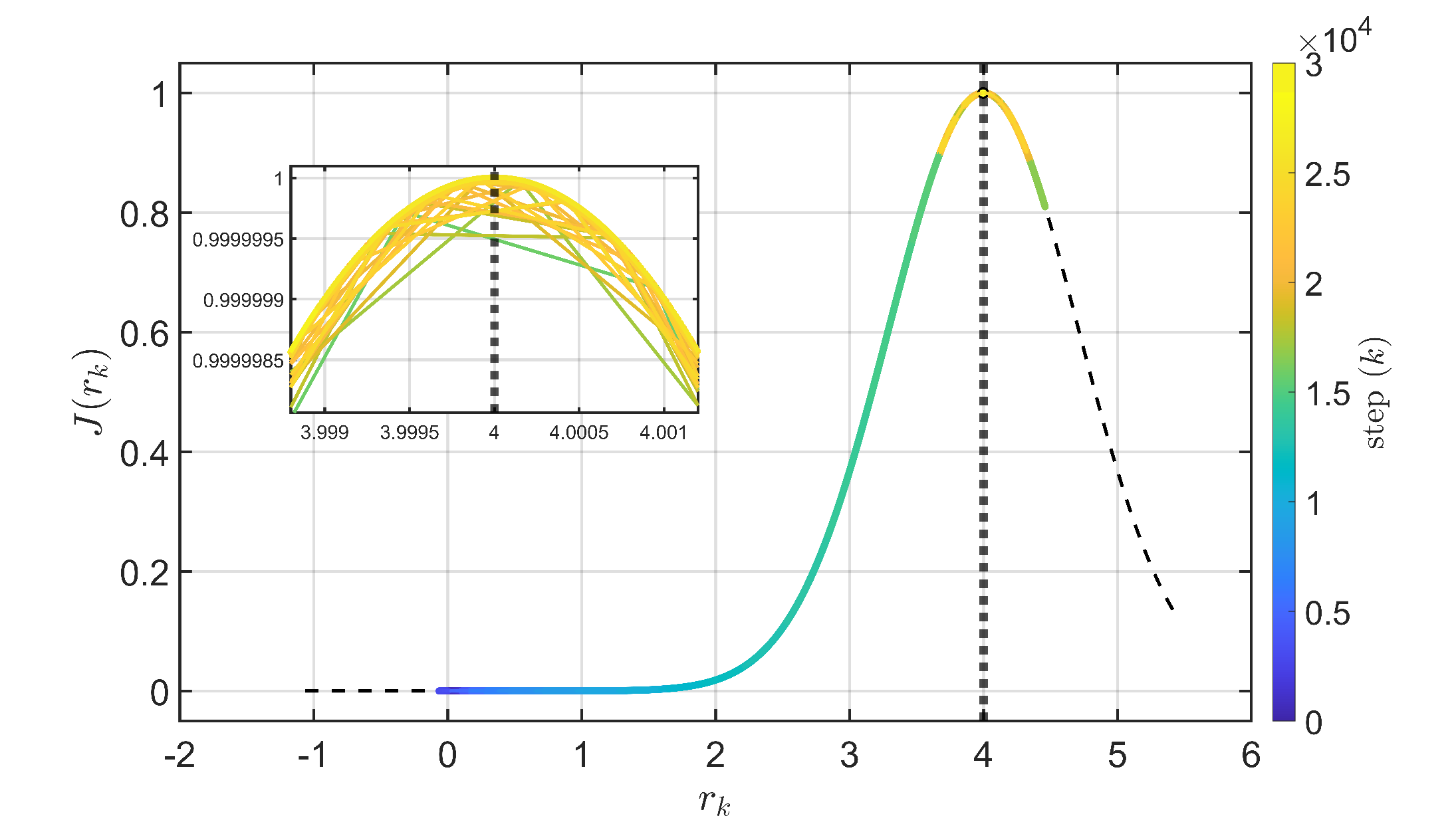}
    \caption{
    \textbf{Exponentially unstable system.} This plot shows the values of $J(r_k)$ computed using the command $r_k$ obtained from \eqref{eq:ECG_output}.
    %
    %
    The step-heat map shows that $r_k$ converges to  $r^\star$ denoted by the vertical dotted line.
    The cost function $J(r_k)$ is shown by the dashed black line.
    }
    \label{fig:Cost_Function_ESC_PCAC_UMS_system}
\end{figure}

The cost function to be maximized is given by
\begin{align}
    J(r)
        &=
            e^
            {
                -(r-r^\star)^2
            },
    \label{eq:USM_cost_function}
\end{align}
which is maximized by $r=r^\star$, where $J(r^\star)=1$.  
For this example, $r^* = 4.$
Figure \ref{fig:Cost_Function_ESC_PCAC_UMS_system} shows that ECG/PCAC updates the command $r$ so that the output of the plant controlled by PCAC converges to $r^*.$  
Although measurements of $J$ are assumed to be available, the functional form of $J$ and the maximizing command $r^*$ are assumed to be unknown.




%

In this example, we use the normalized gradient gain
\begin{align}
    K_{{\rm es},{k}}
        &=
            \dfrac{
                K_{{\rm es},{0}}
            }
            {
                \epsilon
                |y_{{\rm{l}},k}|
            },
    \label{eq:K_esc_normalized}
\end{align}
where
$
    K_{{\rm es},{k}}
$
is the gradient gain at step $k;$
$
    K_{{\rm es},{0}}
$
is the nominal gradient gain; and
$
    \epsilon>0
$
is a normalization constant that prevents division by zero.
In this example, we set
$K_{{\rm es},{0}} = 0.05,$ and $\epsilon = 1e\mbox{-6}.$
The normalization is introduced to address flat regions of the cost function \eqref{eq:USM_cost_function}, where the gradient magnitude can be numerically close to zero not only near the maximizer $r^\star,$ but also when the optimizer is initialized far from $r^\star$ in a region where the cost is nearly flat. 
In such cases, a standard gradient step yields very small updates and may result in slow escape from that region.
The term $\epsilon|y_{{\rm{l}},k}|$ prevents the step size from collapsing when $|y_{{\rm{l}},k}|$ is small, thereby increasing sensitivity to weak gradient information and improving the transient search behavior from poorly initialized conditions.



In this example, the control magnitude-saturation constraint is $u \in [-20,20].$


%
%
%
%

\begin{figure}[!ht]
    \centering
    \includegraphics[width=\columnwidth]{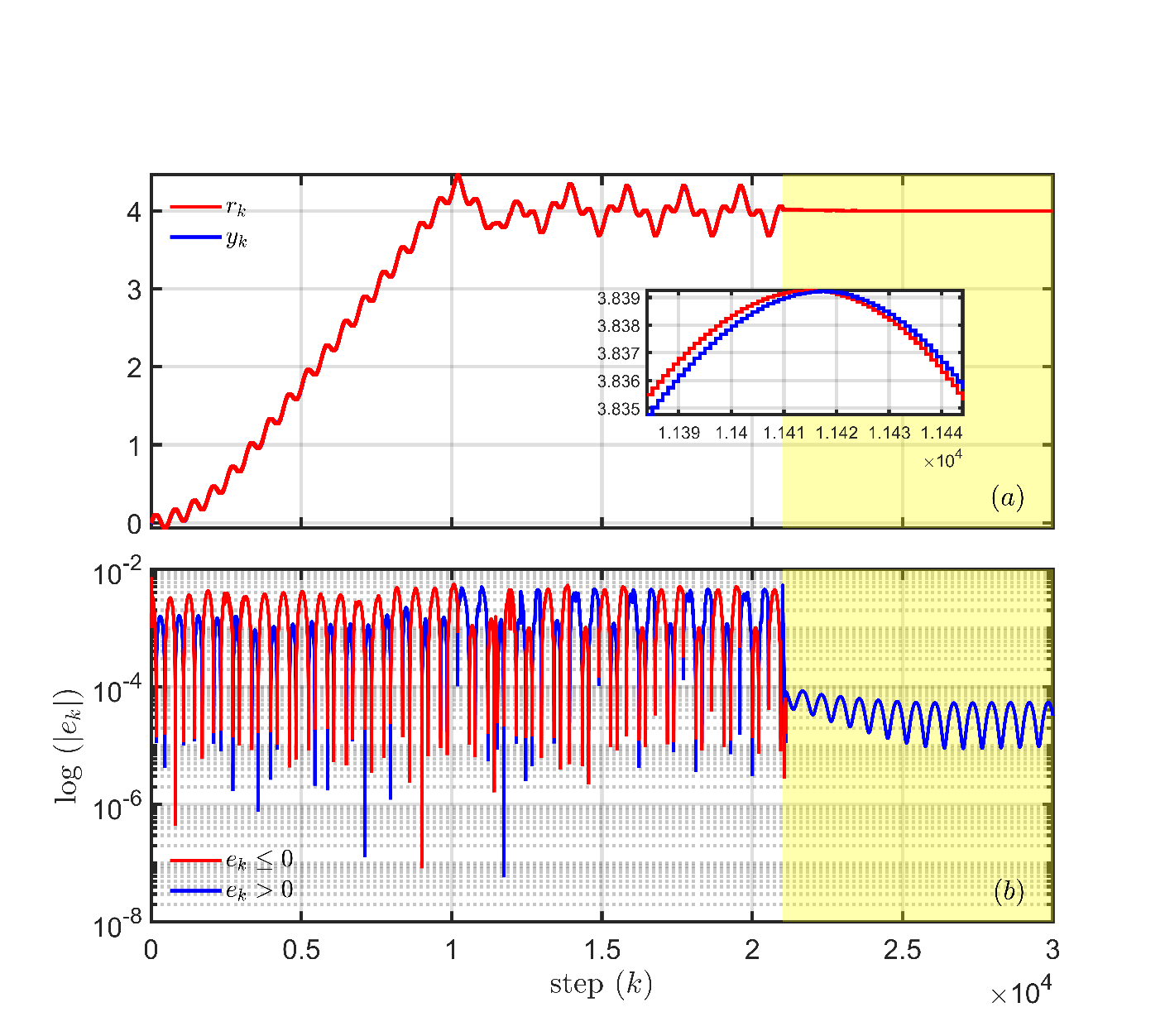}
    \vspace{-.4in}
    \caption{
    \textbf{Exponentially unstable system.} 
    (a) shows the closed-loop response of the exponentially unstable system  \eqref{eq:EUS_x1_dot}, \eqref{eq:EUS_x2_dot} with the control $u_k$ given by PCAC and the command $r_k$ given by ECG.  
    %
    %
    The inset and (b) show that PCAC follows $r_k$.
    %
    }
    \label{fig:ECG_PCAC_Output_UMS_system}
\end{figure}

\begin{figure}[!ht]
    \centering
    \includegraphics[width=\columnwidth]{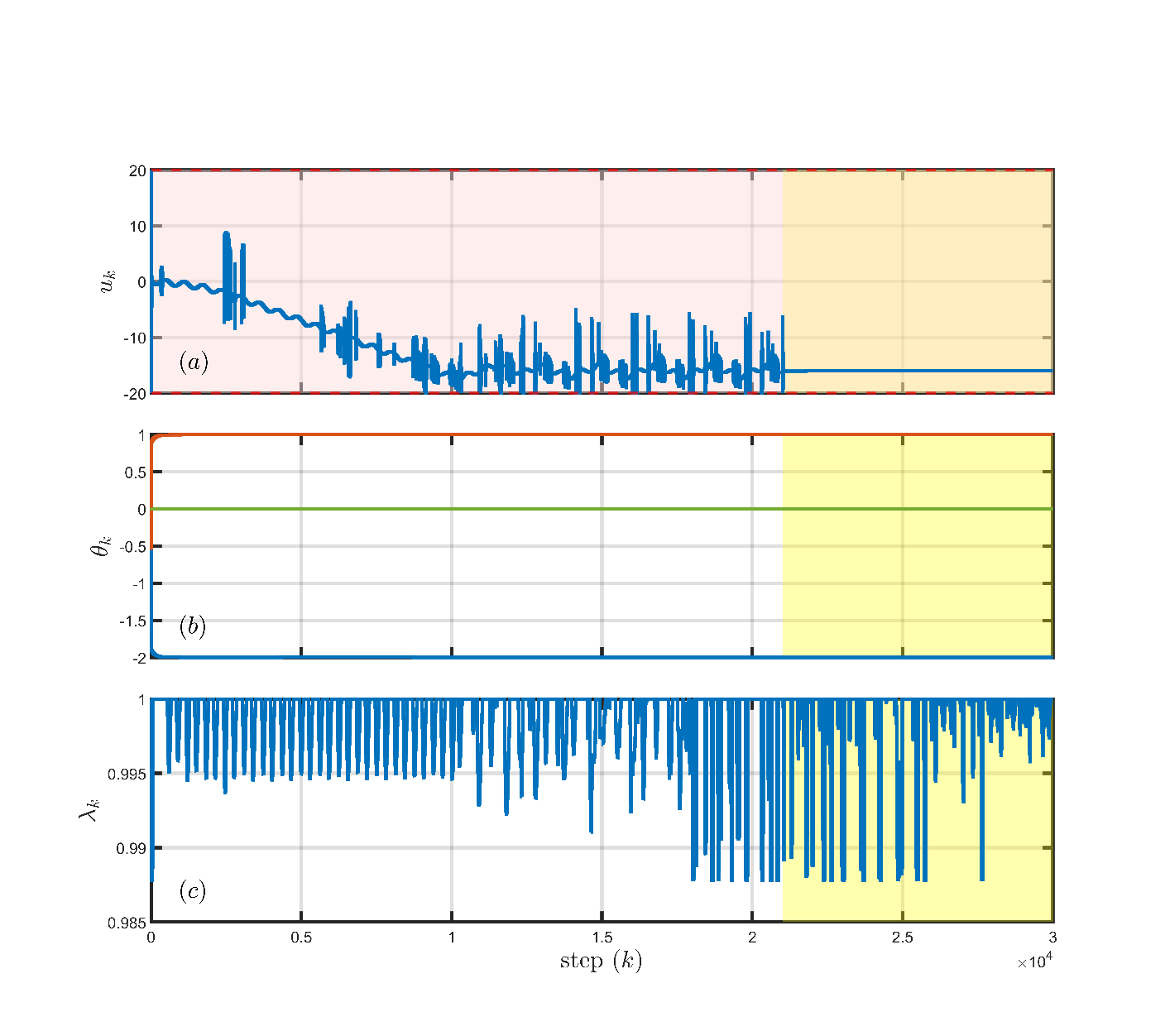}
    \vspace{-.4in}\caption{
    \textbf{Exponentially unstable system.}
    (a) shows the control signal $u_k$ obtained from PCAC,
    (b) shows the corresponding model coefficients obtained from RLS,
    and (c) shows the forgetting factor.
    %
    %
    %
    }
    \label{fig:PCAC_Control_UMS_system}
\end{figure}

For $k > 2.1e4,$ corresponding to the yellow  regions in Figures \ref{fig:ECG_PCAC_Output_UMS_system}, \ref{fig:PCAC_Control_UMS_system}, the dither amplitude and gradient gain are modulated according to
\begin{align}
    a_{{\rm es},{k+1}}
        &=
            \max
            \{
                a_{{\rm es},{\min}}
                , 
                \left(
                    1
                    -
                    \gamma_{\rm a}
                    e^
                    {
                        -
                        \left(
                        \tfrac
                        {
                        y_{{\rm{l}},k}  
                        }
                        {
                        y_{{\rm l},{\rm ref}}
                        }
                        \right)^2
                    }
                \right)
                \,
                a_{{\rm es},{k}}
            \},
    \label{eq:a_es_adaptation_exp_UNS_system}
    \\
    K_{{\rm es},{k+1}}
        &=
            \max
            \{
                K_{{\rm es},{\min}}
                , 
                \left(
                    1
                    -
                    \gamma_{\rm K}
                    e^
                    {
                        -
                        \left(
                        \tfrac
                        {
                        y_{{\rm{l}},k}  
                        }
                        {
                        y_{{\rm l},{\rm ref}}
                        }
                        \right)^2
                    }
                \right)
                \,
                K_{{\rm es},{k}}
            \},
    \label{eq:K_es_adaptation_exp_UNS_system}
\end{align}
where
$a_{{\rm es},{k}}$ and $K_{{\rm es},{k}}$ are the dither amplitude and gain at step $k,$ 
$a_{{\rm es},{\min}} > 0$ and $K_{{\rm es},{\min}} > 0$ are lower bounds that prevent vanishing excitation and gradient gain,
$\gamma_{\rm a}, \gamma_{\rm K} \in (0,1)$ are attenuation gains that regulate the reduction rate,
$y_{{\rm l},{\rm ref}} > 0$ is a normalization constant;
and 
$
    e^
    {
    -
        \left(
            \frac{y_{{\rm l},k}}{y_{{\rm l},{\rm ref}}}
        \right)^2
    }
$ 
provides smooth output-dependent attenuation.
In this example, we set
$a_{{\rm es},{\min}} = K_{{\rm es},{\min}} = 1e\mbox{-3},$   
$\gamma_{\rm a} = \gamma_{\rm K} = 5e\mbox{-2},$ 
and $y_{{\rm l},{\rm ref}} = 0.1.$
These modulations are performed after $2.1e4$ steps to reduce asymptotic oscillations due to the dither signal.
%
%
%
%

Figure \ref{fig:ECG_PCAC_Output_UMS_system} shows the closed-loop response of the exponentially unstable system system \eqref{eq:USM_x1_dot}, \eqref{eq:USM_x2_dot} 
under PCAC, together with the command $r_k$ generated by ECG \eqref{eq:ECG_output} and the command-following error $|e_k|$.
%
%
The yellow region shows that   \eqref{eq:a_es_adaptation_exp_UNS_system}, \eqref{eq:K_es_adaptation_exp_UNS_system}  asymptotically reduces the oscillations.


Figure \ref{fig:PCAC_Control_UMS_system} shows $u_k$ obtained from solving the constrained optimization problem \eqref{cost}--\eqref{lastconstraint},  the estimated model coefficients $\theta_k$ given by \eqref{eq:RLS3}, and the forgetting factor $\lambda_k$ given by \eqref{eq:VRF}.
Figure \ref{fig:Cost_Function_ESC_PCAC_UMS_system} shows the values of $J(r_k)$ computed from \eqref{eq:USM_cost_function} using the command $r_k$ generated by \eqref{eq:ECG_output}, together with its convergence to $r^\star$.







%
%


\section{Conclusions}
\label{sec:conclusions}
This work integrated extremum-seeking-based command generation  (ECG) with predictive cost adaptive control (PCAC) within a closed-loop 
framework that combines command generation with closed-loop system identification and online optimization.
The numerical investigation showed that ESC generates a sequence of commands that approach the optimal command while PCAC ensures stability, command following, and disturbance rejection.
These numerical examples highlight the potential synergy between performance optimization and feedback control without requiring models of the plant dynamics, disturbance, or the performance metric. 
Future work will address extensions to higher order, nonlinear, and MIMO plants.

\printbibliography

\end{document}